\documentclass[10pt]{article}

\usepackage{amsmath}
\usepackage{amsfonts}
\usepackage{amssymb}
\usepackage{framed}
\usepackage{enumerate}
\usepackage{mathrsfs}
\usepackage{enumitem}
\usepackage{comment}
\usepackage[hidelinks]{hyperref}
\usepackage{scrextend}
\usepackage{calligra}
\usepackage{amsthm}
\usepackage{bbm}
\usepackage{thmtools}
\usepackage{thm-restate}
\usepackage{mathabx}
\usepackage[margin=1.5in]{geometry}
\usepackage{rotating}

\usepackage{xcolor}
\hypersetup{
    colorlinks,
    linkcolor={red!50!black},
    citecolor={blue!50!black},
    urlcolor={blue!80!black}
}

\usepackage{tikz-cd}

\theoremstyle{plain}
\newtheorem{thm}{Theorem}[section]
\newtheorem{lem}[thm]{Lemma}
\newtheorem{prop}[thm]{Proposition}
\newtheorem{cor}[thm]{Corollary}

\theoremstyle{definition}
\newtheorem{defn}[thm]{Definition}
\newtheorem{conj}[thm]{Conjecture}

\newtheorem{claim}[thm]{Claim}

\theoremstyle{remark}
\newtheorem*{rem}{Remark}

\newtheorem*{notn}{Notation}

\tikzset{
  symbol/.style={
    draw=none,
    every to/.append style={
      edge node={node [sloped, allow upside down, auto=false]{$#1$}}}
  }
}

\newcommand{\Spec}{\textrm{Spec} \hspace{0.15em} }

\newcommand\restr[2]{{
	\left.\kern-\nulldelimiterspace
	#1
	\vphantom{\big|}
	\right|_{#2}
	}}
\newcommand{\an}[1]{#1^{\textrm{an}}}

\newcommand{\GL}{\textrm{GL}}
\newcommand{\mono}{\textbf{H}}

\tikzset{
    labl/.style={anchor=south, rotate=90, inner sep=.5mm}
}

\DeclareMathOperator{\sheafhom}{\mathscr{H}\text{\kern -3pt {\calligra\large om}}\,}

\title{Absolute Hodge and $\ell$-adic Monodromy}
\author{David Urbanik}

\begin{document}

\maketitle

\begin{abstract}
Let $\mathbb{V}$ be a motivic variation of Hodge structure on a $K$-variety $S$, let $\mathcal{H}$ be the associated $K$-algebraic Hodge bundle, and let $\sigma \in \textrm{Aut}(\mathbb{C}/K)$ be an automorphism. The absolute Hodge conjecture predicts that given a Hodge vector $v \in \mathcal{H}_{\mathbb{C}, s}$ above $s \in S(\mathbb{C})$ which lies inside $\mathbb{V}_{s}$, the conjugate vector $v_{\sigma} \in \mathcal{H}_{\mathbb{C}, s_{\sigma}}$ is Hodge and lies inside $\mathbb{V}_{s_{\sigma}}$. We study this problem in the situation where we have an algebraic subvariety $Z \subset S_{\mathbb{C}}$ containing $s$ whose algebraic monodromy group $\mono_{Z}$ fixes $v$. Using relationships between $\mono_{Z}$ and $\mono_{Z_{\sigma}}$ coming from the theories of complex and $\ell$-adic local systems, we establish a criterion that implies the absolute Hodge conjecture for $v$ subject to a group-theoretic condition on $\mono_{Z}$. We then use our criterion to establish new cases of the absolute Hodge conjecture.
\end{abstract}

\tableofcontents

\section{Introduction}

Let $f : X \to S$ be a smooth projective\footnote{We follow the convention in \cite{Hartshorne1977SchemesChap}, so this means that there exists an embedding $\iota : X \hookrightarrow \mathbb{P}^{n}_{S}$ over $S$.} morphism of $K$-varieties, with $S$ smooth and quasi-projective, and $K \subset \mathbb{C}$ a subfield. In this setting, the local system $\mathbb{V} = R^{2k} f^{\textrm{an}}_{*} \mathbb{Q}(k)$ underlies a polarizable variation of Hodge structure, which has the property that the vector bundle $\mathbb{V} \otimes \mathcal{O}_{\an{S}}$ admits a $K$-algebraic model $\mathcal{H} = R^{2k} f_{*} \Omega^{\bullet}_{X/S}$, where $\Omega^{\bullet}_{X/S}$ is the complex of relative differentials. Moreover, the underlying integral local system $\mathbb{V}_{\mathbb{Z}} \subset \mathbb{V}$ admits a comparison with the $\ell$-adic local system $\mathbb{V}_{\ell} = R^{2k} f^{\textrm{\'et}}_{*} \mathbb{Z}_{\ell}(k)$ on the \'etale site of $S_{\mathbb{C}}$. 

\smallskip

\begin{notn}
For $s \in S(\mathbb{C})$, we denote by $c^{\textrm{dR}}_{s} : \mathbb{V}_{\mathbb{C}, s} \xrightarrow{\sim} \mathcal{H}_{\mathbb{C}, s}$ and $c^{\ell}_{s} : \mathbb{V}_{\mathbb{Q}_{\ell}, s} \xrightarrow{\sim} \mathbb{V}_{\ell, \mathbb{Q}_{\ell}, s}$ the natural comparison isomorphisms.
\end{notn}


\smallskip

The fact that $\mathcal{H}$ and $\mathbb{V}_{\ell}$ come from a morphism $f : X \to S$ of $K$-varieties means that for each automorphism $\sigma \in \textrm{Aut}(\mathbb{C}/K)$ and point $s \in S(\mathbb{C})$, there are isomorphisms $(-)_{\sigma} : \mathcal{H}_{\mathbb{C}, s} \xrightarrow{\sim} \mathcal{H}_{\mathbb{C}, s_{\sigma}}$ and $(-)_{\sigma} : \mathbb{V}_{\ell, \mathbb{Q}_{\ell}, s} \xrightarrow{\sim} \mathbb{V}_{\ell, \mathbb{Q}_{\ell}, s_{\sigma}}$ induced by conjugation. Combining these isomorphisms with the isomorphisms $c^{\textrm{dR}}_{s}$ and $c^{\ell}_{s}$, we may consider several possible ways in which the rational structure of $\mathbb{V}$ may be preserved under conjugation:

\smallskip

\begin{defn}
\label{conjdef}
Let $s \in S(\mathbb{C})$ be a point, $v \in \mathbb{V}_{s}$ be a vector, $W \subset \mathbb{V}_{s}$ a subspace, and $\sigma \in \textrm{Aut}(\mathbb{C}/K)$. Then we say:
\begin{itemize}
\item[(i)] That $v$ has \emph{rational conjugates} if $c_{s}^{\textrm{dR}}(v)_{\sigma}$ (resp. $c_{s}^{\ell}(v)_{\sigma}$ for all $\ell$) are rational vectors of $\mathbb{V}$. (They lie in the image of $\mathbb{V}_{s_{\sigma}}$ under the comparison isomorphisms $c_{s_{\sigma}}^{\textrm{dR}}$ and $c^{\ell}_{s_{\sigma}}$.)
\item[(ii)] That $v$ has a \emph{canonical rational conjugate} if there exists $v_{\sigma} \in \mathbb{V}_{s_{\sigma}}$ such that $c_{s}^{\textrm{dR}}(v)_{\sigma} = c_{s_{\sigma}}^{\textrm{dR}}(v_{\sigma})$ and $c_{s}^{\ell}(v)_{\sigma} = c_{s_{\sigma}}^{\ell}(v_{\sigma})$ for all $\ell$. (The vector $v_{\sigma}$ works for all comparison isomorphisms at once.)
\item[(iii)] That $W$ has \emph{rational conjugates} if each element of $W$ does.
\item[(iv)] That $W$ has a \emph{canonical rational conjugate} if each element of $W$ does.
\end{itemize}
\end{defn}

\begin{rem}
It is clear that if $w$ and $w'$ have (canonical) rational conjugates then so do $w + w'$ and $\lambda w$ for any $\lambda \in \mathbb{Q}$, so we lose nothing by considering subspaces with (canonical) rational conjugates.
\end{rem}

\smallskip

In the case $K = \mathbb{Q}$, the property that a Hodge vector $v$ has rational conjugates is equivalent to the absolute Hodge conjecture for $v$ formulated in \cite{ChapterElevenNotesonAbsoluteHodgeClasses}.\footnote{Although the formulation in \cite{Deligne1982}, notably, is equivalent to the condition that $v$ has a \emph{canonical} rational conjugate; see also the discussion in \autoref{abshodgesec}.} The property that $v$ has a canonical rational conjugate is a priori stronger, and would follow for Hodge vectors from the Hodge conjecture in the geometric setting. These properties are expected to hold for Hodge vectors more generally in the \emph{motivic} setting. Loosely speaking, we will say that a variation $\mathbb{V}$ is motivic if it is generated from the geometric situation by functorial constructions; we give a precise definition suitable for our purposes in \autoref{propertiessec}.

Our central technical result, from which new cases of the absolute Hodge conjecture will follow, gives a new way of deducing that a subspace $W \subset \mathbb{V}_{s}$ (or a subspace of the various tensor powers of $\mathbb{V}_{s}$) has rational conjugates by studying the monodromy of complex subvarieties $Z \subset S_{\mathbb{C}}$. To formulate the statement we introduce some notation.

\smallskip

\begin{notn}
Given an $R$-module $V$, we let $V^{m,n} = (V^{*})^{\otimes m} \otimes V^{\otimes n}$, with the analogous definition for $R$-local systems $\mathbb{V}$.
\end{notn}

\begin{notn}
Given an $R$-module $V$, we let $\mathcal{T}(V) = \bigoplus_{m, n \geq 0} V^{n,m}$, with the analogous definition for $R$-local systems $\mathbb{V}$.
\end{notn}

\begin{notn}
Let $V$ be a vector space over $R$, and let $T \subset \mathcal{T}(V)$ be a vector subspace. We denote by $\mathbf{G}_{T} \subset \textrm{Aut}(V)$ the $R$-algebraic subgroup defined by the property that it fixes each element of $T$.
\end{notn}

\begin{notn}
Let $Z \subset S_{\mathbb{C}}$ be an algebraic subvariety, and $\mathbb{V}$ be a local system. We denote by $\mono_{Z}$ the \emph{algebraic monodromy} of $Z$, which is the identity component of the Zariski closure of the monodromy representation on $\restr{\mathbb{V}}{Z^{\textrm{nor}}}$, where $Z^{\textrm{nor}} \to Z$ is the normalization of $Z$.
\end{notn}

\begin{rem}
Given a subvariety $Z \subset S_{\mathbb{C}}$, we may regard $\mono_{Z}$ as an algebraic subgroup of $\textrm{Aut}(\mathbb{V}_{s})$ for any point $s \in Z(\mathbb{C})$. We will use the notation $\mono_{Z}$ for each of these subgroups, even though the point $s$ may vary, as the context will leave no ambiguity about which algebraic group we mean.
\end{rem}

\smallskip

Our main technical result is then the following:

\begin{defn}
If $G \subset \textrm{Aut}(\mathbb{V}_{s})$ is an algebraic group, and $W \subset \mathcal{T}(\mathbb{V}_{s})$ is a subspace, we will say that $G$ fixes $W$ as a $\mathbb{Q}$-subspace if for each $g \in G(\mathbb{C})$ we have $g W = W$.\footnote{We wish to stress that this condition is different from fixing $W$ pointwise, or fixing $W_{\mathbb{C}}$ as a subspace.}
\end{defn}

\begin{thm}
\label{mainthm}
Let $\mathbb{V}$ be a motivic variation of Hodge structure on $S$, and let $s \in S(\mathbb{C})$. Suppose that $T \subset \mathcal{T}(\mathbb{V}_{s})$ is a subspace with canonical rational conjugates, that $Z \subset S_{\mathbb{C}}$ is a complex subvariety containing $s$ satisfying $\mono_{Z} \subset \mathbf{G}_{T}$, and that $W \subset \mathcal{T}(\mathbb{V}_{s})$ is a subspace containing $T$. Suppose that the normalizer $N$ of $\mono_{Z}$ inside $\mathbf{G}_{T}$ fixes $W$ as a $\mathbb{Q}$-subspace. Then $W$ has rational conjugates.
\end{thm}

A simple application of \autoref{mainthm} is then: 

\begin{cor}
\label{maincor}
Suppose $\mathbb{V}$ is a motivic variation of Hodge structure on $S$, let $s \in S(\mathbb{C})$ and suppose that $Z \subset S_{\mathbb{C}}$ contains $s$. Then if the only $(m,n)$ tensors fixed by $\mono_{Z}$ are of the form $\mathbb{Q} w$ for some $w \in \mathbb{V}_{s}^{m,n}$, then the subspace $\mathbb{Q} w$ has rational conjugates.
\end{cor}

\begin{proof}
We show in \autoref{polconstr} that any motivic variation admits a polarization $Q : \mathbb{V} \otimes \mathbb{V} \to \mathbb{Q}$ such that the class of $Q_{s}$ inside $\mathbb{V}^{2,0}_{s}$ has a canonical rational conjugate. Take $T$ to be the span of $Q_{s}$. From the fact that $N$ preserves $Q$ we learn it acts on $w$ by $\pm 1$. Taking $W$ to be the span of $Q_{s}$ and $w$, we learn that $W$ has rational conjugates.
\end{proof}

\noindent As we explain in \autoref{voisinexamplesec}, a form of the argument in \autoref{maincor} appears in Voisin's paper \cite{voisin_2007}. There it is used to prove a ``weakly absolute Hodge'' statement for Hodge classes coming from a family of fourfolds for which neither the Hodge nor absolute Hodge conjecture are known. Our \autoref{mainthm} gives the full absolute conjecture in this case; this is \autoref{voisinexample}.

\vspace{0.5em}

Let $\textrm{Isom}(\mathbb{V}_{s}, \mathbb{V}_{s_{\sigma}})$ be the $\mathbb{Q}$-algebraic variety consisting of all isomorphisms between $\mathbb{V}_{s}$ and $\mathbb{V}_{s_{\sigma}}$. Our proof of \autoref{mainthm} proceeds through the study of the $\mathbb{Q}$-subvariety
\[ \mathcal{I}_{T}(s, Z, \sigma) \subset \textrm{Isom}(\mathbb{V}_{s}, \mathbb{V}_{s_{\sigma}}) , \]
defined by the property that $r \in \mathcal{I}_{T}(s, Z, \sigma)(\mathbb{C})$ if and only if $r \circ \mono_{Z, \mathbb{C}} \circ r^{-1} = \mono_{Z_{\sigma}, \mathbb{C}}$ and $r$ sends each $t \in T$ to its canonical rational conjugate $t_{\sigma} \in T_{\sigma}$. We observe that $\mathcal{I}_{T}(s, Z, \sigma)$ is a torsor under the normalizer $N$ of $\mono_{Z}$ inside $\mathbf{G}_{T}$. \autoref{mainthm} follows immediately from the following two facts about $\mathcal{I}_{T}(s, Z, \sigma)$, established in \autoref{pointsofIsec}:

\newpage

\begin{prop}
\label{torsorpoints}
~
\begin{itemize}
\item[(i)] The map
\[ r^{\textrm{dR}}_{\sigma} : \mathbb{V}_{\mathbb{C}, s} \xrightarrow{c^{\textrm{dR}}_{s}} \mathcal{H}_{\mathbb{C}, s} \xrightarrow{\sigma} \mathcal{H}_{\mathbb{C}, s_{\sigma}} \xrightarrow{(c^{\textrm{dR}}_{s_{\sigma}})^{-1}} \mathbb{V}_{\mathbb{C}, s_{\sigma}} \]
induced by $\sigma$ is a $\mathbb{C}$-point of $\mathcal{I}_{T}(s, Z, \sigma)$. 
\item[(ii)] For each $\ell$, the map 
\[ r^{\ell}_{\sigma} : \mathbb{V}_{\mathbb{Q}_{\ell}, s} \xrightarrow{c^{\ell}_{s}} \mathbb{V}_{\ell, \mathbb{Q}_{\ell}, s} \xrightarrow{\sigma} \mathbb{V}_{\ell, \mathbb{Q}_{\ell}, s_{\sigma}} \xrightarrow{(c^{\ell}_{s_{\sigma}})^{-1}} \mathbb{V}_{\mathbb{Q}_{\ell}, s_{\sigma}} \]
induced by $\sigma$ is a $\mathbb{Q}_{\ell}$-point of $\mathcal{I}_{T}(s, Z, \sigma)$.
\end{itemize}
\end{prop}

\vspace{0.5em}

\begin{proof}[Proof of \autoref{mainthm} (assuming \ref{torsorpoints}):] 

~

\vspace{0.5em}

We let $\mathcal{I} = \mathcal{I}_{T}(s, Z, \sigma)$ for ease of notation. That $\mathcal{I}$ is naturally a torsor under $N$ implies that the subspace $r(W)$ is independent of the choice of $r \in \mathcal{I}(\mathbb{C})$: any two choices $r$ and $r'$ are related by $r' = r n$ with $n \in N(\mathbb{C})$, and $N$ fixes $W$ by assumption. As $\mathcal{I}$ is a $\mathbb{Q}$-variety, it has a $\overline{\mathbb{Q}}$-point $r$, hence $W_{\sigma} := r(W)$ is defined over $\overline{\mathbb{Q}}$. If we can show that $W_{\sigma}$ is in fact a $\mathbb{Q}$-subspace of $\mathbb{V}_{s_{\sigma}}$, then the proof will be complete as \autoref{torsorpoints} then implies both that $c^{\textrm{dR}}_{s}(W)_{\sigma} = c^{\textrm{dR}}_{s_{\sigma}}(W_{\sigma})$ and $c^{\ell}_{s}(W)_{\sigma} = c^{\ell}_{s_{\sigma}}(W_{\sigma})$ for each $\ell$.

By \autoref{torsorpoints}(ii) the variety $\mathcal{I}$ has $\mathbb{Q}_{\ell}$-points for every $\ell$, hence it follows that under some choice of isomorphism $\mathbb{C} \cong \overline{\mathbb{Q}_{\ell}}$, every element of $W_{\sigma}$ is defined over $\mathbb{Q}_{\ell}$. The proof is completed by the following:

\begin{claim}
\label{affinespptclaim}
Let $v \in \mathbb{A}^{n}(\overline{\mathbb{Q}})$ be a $\overline{\mathbb{Q}}$-point in affine space, and suppose that for each $\ell$ there exists an isomorphism $\mathbb{C} \cong \overline{\mathbb{Q}_{\ell}}$ under which $v$ is identified with a $\mathbb{Q}_{\ell}$-point. Then $v$ is defined over $\mathbb{Q}$.
\end{claim}
\begin{proof}
As $\mathbb{A}^{n} = (\mathbb{A}^{1})^{n}$, it suffices to consider the case $n = 1$. Let $f$ be the minimal polynomial of $v$. The statement is unchanged by replacing $v$ with $\lambda v$ for $\lambda \in \mathbb{Q}$, so we may assume that $v$ is integral and that $f$ is monic with integral coefficients. The hypothesis shows (using integrality) that $f$ has a root over $\mathbb{\mathbb{Z}}_{\ell}$ for every $\ell$, hence a root modulo $\ell$ for every $\ell$. As $f$ is irreducible, this implies that $f$ is linear by \cite[p.~362,~Ex.~6.2]{cassels1967algebraic}.
\end{proof}
\end{proof}

\begin{rem}
The proof of \autoref{mainthm} does not invoke the properties of the filtration $F^{\bullet}$ in any way, so it applies to motivic mixed variations as well. We focus exclusively on the pure case in this paper.
\end{rem}

\vspace{1em}

\autoref{mainthm} lets us deduce the following result on the absolute Hodge conjecture, as we show in \autoref{abshodgesec}: 

\begin{thm}
\label{mainthm2}
Let $\mathbb{V}$, $S$, $s$, $T$, and $Z$ be as in \autoref{mainthm}, and suppose that $W$ is a subspace of Hodge vectors. Then the $K$-absolute Hodge conjecture\footnote{See \autoref{abshodgesec} for a definition; the usual absolute Hodge conjecture is the case $K = \mathbb{Q}$.} holds for each element of $W$.
\end{thm}

We give two applications of our results. The first is to resolve the absolute Hodge conjecture for classes coming from a family of fourfolds considered by Voisin in her paper \cite[Example~3.4]{voisin_2007}. More specifically, we show the following:
\begin{thm}
\label{voisinexample}
Let $f : X \to S$ be the family of degree $6$ hypersurfaces in $\mathbb{P}^{5}$ stable under the involution $\iota(X_{0}, \hdots, X_{5}) = (-X_{0}, -X_{1}, X_{2}, \hdots, X_{5})$, let $\mathbb{V}$ be the geometric variation on middle cohomology of degree $2k$, and let $\mathbb{V}_{-}$ be the subvariation obtained by taking anti-invariants under $\iota$. Then if $\mathcal{H}_{-}^{k,k}$ is the quotient $F_{-}^{k}/F_{-}^{k+1}$, the rational Hodge classes $v \in \mathbb{V}_{-}$ for which the absolute Hodge conjecture holds are topologically dense in the real part of $\mathcal{H}_{-}^{k,k}$.
\end{thm}
\noindent Our method of establishing \autoref{voisinexample} is similar to the argument employed in Voisin's paper to establish that such classes are ``weakly absolute Hodge''. We are able to obtain a stronger conclusion due to the presence of $\ell$-adic input in our argument. 

As a second application, we let $\mathbb{V}$ be a motivic variation on $S$, let $Z \subset S \times S$ be the algebraic locus of points $(x, y)$ such that $\mathbb{V}_{x} \cong \mathbb{V}_{y}$ as integral polarized Hodge structures, and we let $Z_{0} \subset Z$ be an irreducible component containing the diagonal $\Delta \subset S \times S$. The variety $Z_{0}$ is a component of the Hodge locus for the variation $\mathbb{V} \boxtimes \mathbb{V}^{*}$ on $S \times S$ (see \autoref{boxvardef}). We call an isomorphism $\mathbb{V}_{x} \cong \mathbb{V}_{y}$ with $(x,y) \in Z_{0}$ \emph{generic} if its class is Hodge at a generic point of $Z_{0}$. We have the following result:


\begin{thm}
\label{genisothm}
If $\mathbb{V}$ is a motivic variation of Hodge structure on the smooth $K$-variety $S$ such that $\mono_{S} = \textrm{Aut}(\mathbb{V}_{s}, Q_{s})$ (monodromy is Zariski dense), then rational tensors fixed by $\mono_{Z_{0}}$ have rational conjugates, and hence generic isomorphisms between $\mathbb{V}_{x}$ and $\mathbb{V}_{y}$ are $K$-absolutely Hodge.
\end{thm}

\subsection{Structure of the Paper}

In \autoref{motivicvarsec} we discuss motivic variations and the Hodge and absolute Hodge conjectures for such variations. We review some properties of motivic variations in \autoref{propertiessec}, and list all the properties that we will be needed for the remainder of the paper. The properties in \autoref{propertiessec} relating to polarizations require further justification, so we carry out the necessary constructions in  \autoref{polconstr}. As mentioned, \autoref{abshodgesec} discusses the relationship between our notion of rational conjugates and the Absolute Hodge conjecture. 

In \autoref{pointsofIsec} we prove \autoref{torsorpoints}. The proof of \autoref{torsorpoints}(i) is already implicit in the paper \cite{KOU}, so we give a brief summary of the argument in \autoref{derhamconjproof}. To tackle the $\ell$-adic case in \autoref{ladicconjproof} we reduce the necessary statement to a purely algebraic statement involving \'etale fundamental groups, after which the result follows by showing that $\sigma \in \textrm{Aut}(\mathbb{C}/K)$ essentially conjugates the representation of $\pi_{1}^{\textrm{\'et}}(S_{\mathbb{C}}, s)$ on the fibre $\mathbb{V}_{\ell, s}$ to the representation of $\pi_{1}^{\textrm{\'et}}(S_{\mathbb{C}}, s_{\sigma})$ on the fibre $\mathbb{V}_{\ell, s_{\sigma}}$.

The same ideas that give \autoref{mainthm} and \autoref{mainthm2} in fact give several other variants of \autoref{mainthm} and \autoref{mainthm2}, as we discuss in \autoref{mainthmsec}. We then consider the example of Voisin in \autoref{voisinexamplesec}, and the application to generic isomorphisms of Hodge structures in \autoref{genisosec}.

\section{Motivic Variations}
\label{motivicvarsec}

\subsection{Properties of Motivic Variations}
\label{propertiessec}

Recall that our central example of a variation of Hodge structures is the variation $\mathbb{V}$ obtained from a smooth, projective family $f : X \to S$ of $K$-varieties, with $S$ smooth. Let us briefly state some properties of such variations:

\begin{itemize}
\item[(i)] The Hodge bundle $\mathbb{V} \otimes \mathcal{O}_{\an{S}}$ admits a canonical model as a $K$-algebraic vector bundle $\mathcal{H} = R^{2k} f_{*} \Omega^{\bullet}_{X/S}$, where $\Omega^{\bullet}_{X/S}$ is the complex of relative algebraic differentials.
\item[(ii)] The filtration $F^{\bullet}$ on $\mathcal{H}$ admits a canonical $K$-algebraic model, agreeing with the filtration in the analytic setting. This is obtained from a filtration on $\Omega^{\bullet}_{X/S}$ and an appropriate spectral sequence.
\item[(iii)] The connection $\nabla$ whose flat sections are given by $\mathbb{V}_{\mathbb{C}}$ is $K$-algebraic. This is due to Katz and Oda \cite{katz1968}.
\item[(iv)] The $\mathbb{Q}$-local system $\mathbb{V}$ admits a canonical integral subsystem $\mathbb{V}_{\mathbb{Z}} \subset \mathbb{V}$, such that $\mathbb{V}_{\mathbb{Z}} \otimes \mathbb{Z}_{\ell}$ admits a comparison isomorphism with $\an{\mathbb{V}_{\ell}}$, where $\mathbb{V}_{\ell}$ is an $\ell$-adic local system defined on the \'etale site of $S_{\mathbb{C}}$. This is the relative version of the comparison between Betti and \'etale cohomology.
\end{itemize}

\noindent The following additional properties relate to the polarization on $\mathbb{V}$, and will be established in the next section:

\begin{itemize}
\item[(v)] The variation $\mathbb{V}$ admits a polarization $Q : \mathbb{V} \otimes \mathbb{V} \to \mathbb{Q}$ which is $K$-algebraic, in the sense that we have a morphism $\mathcal{Q} : \mathcal{H} \otimes \mathcal{H} \to \mathcal{O}_{S}$ of $K$-algebraic vector bundles whose analytification is identified with $Q$ after applying the isomorphism $\mathbb{V} \otimes \mathcal{O}_{\an{S}} \simeq \an{\mathcal{H}}$.
\item[(vi)] There exists a bilinear form $Q_{\ell} : \mathbb{V}_{\ell} \otimes \mathbb{V}_{\ell} \to \mathbb{Z}_{\ell}$, compatible with the polarization $Q$ in (v).
\item[(vii)] Let $\sigma \in \textrm{Aut}(\mathbb{C}/K)$, and denote by $\tau_{\sigma} : S_{\textrm{\'et}} \to S_{\textrm{\'et}}$ the induced automorphism of the \'etale site of $S_{\mathbb{C}}$. Then there exists a canonical isomorphism $\tau_{\sigma}^{*} \mathbb{V}_{\ell} \simeq \mathbb{V}_{\ell}$, compatible with the polarization $Q_{\ell}$ in (vi).
\end{itemize}

The properties (i) through (vii) are preserved under any sufficiently functorial construction; for instance, the properties (i) through (vii) are preserved under duals, direct sums and tensor products. The theorems we prove will only ever use the above listed properties, and we will refer to such variations as \emph{motivic}. Note that it is not immediately clear that this notion of motivic variation is the same as other similar notions that appear in the literature, but it is the definition which will be most useful for us.

We note the following elementary consequence of (v), (vi) and (vii) which will be useful later: 

\begin{lem}
\label{conjpreservespol}
For motivic variations, the $\sigma$-conjugate of the class $Q_{s} \in \mathbb{V}^{2,0}_{s}$ is rational for $\mathbb{V}^{2,0}_{s_{\sigma}}$ (and given by $Q_{s_{\sigma}}$) in both the de Rham and $\ell$-adic setting, where $\sigma \in \textrm{Aut}(\mathbb{C}/K)$.
\qed
\end{lem}

\subsection{Construction of $K$-algebraic Polarizations}
\label{polconstr}

The material in this section should be well-known to experts, but we cannot find an appropriate reference. 

We consider the geometric situation with $\mathbb{V} = R^{2k} f_{\textrm{an}, *} \mathbb{Z}$, $\mathcal{H} = R^{2k} f_{*} \Omega^{\bullet}_{X/S}$ and $\mathbb{V}_{\ell} = R^{2k} f_{\textrm{\'et}, *} \mathbb{Z}_{\ell}$. We will define the primitive cohomology subsystems $\mathbb{V}_{\textrm{prim}} \subset \mathbb{V}$ and $\mathbb{V}_{\ell, \textrm{prim}} \subset \mathbb{V}_{\ell}$, and vector subbundle $\mathcal{H}_{\textrm{prim}} \subset \mathcal{H}$, as well as polarizations $Q, Q_{\ell}$ and $\mathcal{Q}$ on the primitive cohomology satisfying the properties (v), (vi) and (vii). That properties (v), (vi) and (vii) hold for the original local systems and vector bundles follows from the usual procedure of obtaining a polarization on all of cohomology from a polarization on the primitive part.

In the case where $S = \Spec K$ is a point, our definition of $Q$ and $\mathcal{Q}$ essentially appears (among other places) in \cite{hodgeII}. Our goal is to generalize this to the relative setting, where we regard the scheme $X$ as an $S$-scheme via the map $f$. The first task is to define the relative Chern classes $c^{\textrm{dR}}_{1, S}(\mathcal{L})$ and $c^{\textrm{\'et}}_{1, S}(\mathcal{L})$ of a line bundle $\mathcal{L}$ on $X$.

\begin{lem}
\label{globaldRchern}
Let $\mathcal{L}$ be a line bundle on $X$. Then there exists a (necessarily unique) global section $c^{\textrm{dR}}_{1, S}(\mathcal{L})$ of $R^{2} f_{*} \Omega^{\bullet}_{X/S}$ whose restriction to the fibre at $s \in S$ is the Chern class $c^{\textrm{dR}}_{1}\left(\restr{\mathcal{L}}{X_{s}}\right)$.\footnote{In the sense that the equality holds after identifying $(R^{2} f_{*} \Omega^{\bullet}_{X/S})_{s}$ with the algebraic de Rham cohomology of the fibre $X_{s}$ using proper base change.}
\end{lem}

\begin{proof}
By \cite[\href{https://stacks.math.columbia.edu/tag/0FLE}{Section 0FLE}]{stacks-project} we obtain a class $c_{1}^{\textrm{dR}}(\mathcal{L}) \in H^{2}_{\textrm{dR}}(X/S)$, where $H^{\bullet}_{\textrm{dR}}(X/S)$ is the cohomology of the complex $R \, \Gamma \, \Omega^{\bullet}_{X/S}$. Choose an injective resolution $\Omega^{\bullet}_{X/S} \to \mathcal{I}^{\bullet}$. We may compute $R^{2} f_{*} \Omega^{\bullet}_{X/S}$ as the cohomology sheaf of $f_{*} \mathcal{I}^{\bullet}$ in degree two. From the equality $\Gamma(S, f_{*} \mathcal{I}^{\bullet}) = \Gamma(X, \mathcal{I}^{\bullet})$ we therefore obtain a morphism $g : H^{2}_{\textrm{dR}}(X/S) \to \Gamma(S, R^{2} f_{*} \Omega^{\bullet}_{X/S})$, and we may define $c^{\textrm{dR}}_{1, S}(\mathcal{L})$ as the image of $c_{1}^{\textrm{dR}}(\mathcal{L})$ under $g$. 

Let $s \in S$ be a point, and $r : \Gamma(S, R^{2} f_{*} \Omega^{\bullet}_{X/S}) \to (R^{2} f_{*} \Omega^{\bullet}_{X/S})_{s}$ be the restriction of a global section to the fibre at $s$. Let $b : (R^{2} f_{*} \Omega^{\bullet}_{X/S})_{s} \to H^{2}_{\textrm{dR}}(X_{s}/\kappa(s))$ be the canonical base change morphism. As Chern classes are functorial, the map $h : H^{2}_{\textrm{dR}}(X/S) \to H^{2}_{\textrm{dR}}(X_{s}/\kappa(s))$ induced by the inclusion $\iota: X_{s} \hookrightarrow X$ sends $c^{dR}_{1}(\mathcal{L})$ to $c^{\textrm{dR}}_{1}\left(\restr{\mathcal{L}}{X_{s}}\right)$. It therefore suffices to check that $h = b \circ r \circ g$.

To make this verification, we recall the construction of the base change map at the level of complexes, following \cite[\href{https://stacks.math.columbia.edu/tag/0735}{Section 0735}]{stacks-project}. We may identify $\Omega^{\bullet}_{X_{s}}$ with the pullback $\iota^{*} \Omega^{\bullet}_{X/S}$, and find an injective resolution $\iota^{*} \Omega^{\bullet}_{X/S} \to \mathcal{J}^{\bullet}$. We then obtain a commuting diagram
\begin{center}
\begin{tikzcd}
\iota_{*} \iota^{*} \Omega^{\bullet}_{X/S} \arrow[r] & 
\iota_{*} \mathcal{J}^{\bullet} \\
\Omega^{\bullet}_{X/S} \arrow[u, "\textrm{adj.}"] \arrow[r] & \mathcal{I}^{\bullet} \arrow[u, "\beta"] 
\end{tikzcd} ,
\end{center}
where the map $\beta$ is unique up to homotopy. The map $b$ (resp. the map $h$) is constructed from $\beta$ by applying $f_{*}$ (resp. applying $\Gamma$) and taking the induced map on cohomology. The required equality then follows from the fact that $(f \circ \iota)_{*} = \Gamma$.
\end{proof}

\begin{lem}
\label{globaletchern}
Let $\overline{K}$ be an algebraic closure of $K$ inside $\mathbb{C}$, and let $\mathcal{L}$ be a line bundle on $X_{\overline{K}}$. Then there exists a (necessarily unique) global section $c^{\textrm{\'et}}_{1, S_{\overline{K}}}(\mathcal{L})$ of $R^{2} f_{\textrm{\'et}, *} \mathbb{Z}_{\ell}(1)$ whose restriction to the fibre at $s \in S(\overline{K})$ is the Chern class $c^{\textrm{\'et}}_{1}\left(\restr{\mathcal{L}}{X_{s}}\right)$.\footnote{In the same sense as in \autoref{globaldRchern}.}
\end{lem}

\begin{proof}
We argue analogously to \autoref{globaldRchern}. For ease of notation, we assume $K = \overline{K}$ and so replace $X_{\overline{K}}$ by $X$ and $S_{\overline{K}}$ by $S$. By \cite[VI, \S 10]{milneetcoh} we obtain a class $c^{\textrm{\'et}}_{1}(\mathcal{L}) \in H^{2}_{\textrm{\'et}}(X, \mathbb{Z}_{\ell}(1))$. Choose an injective resolution $\mathbb{Z}_{\ell} \to \mathcal{I}^{\bullet}$.\footnote{Here we really mean a compatible system of resolutions $\mathbb{Z}/\ell^{k}\mathbb{Z} \to \mathcal{I}_{k}$, where we regard $\ell$-adic sheaves as systems of $\mathbb{Z}/\ell^{k}\mathbb{Z}$-sheaves in the usual way.} Proceeding as before, we may compute $H^{2}_{\textrm{\'et}}(X, \mathbb{Z}_{\ell}(1))$ from the degree two cohomology of $\mathcal{I}^{\bullet}$, and $\mathbb{V}^{2}_{\ell}$ as the degree two cohomology sheaf of $f_{\textrm{\'et},*} \mathcal{I}^{\bullet}$. From the equality $\Gamma(S, f_{\textrm{\'et},*} \mathcal{I}^{\bullet}) = \Gamma(X, \mathcal{I}^{\bullet})$ we therefore obtain a map $g : H^{2}_{\textrm{\'et}}(X, \mathbb{Z}_{\ell}(1)) \to\Gamma(S, R^{2} f_{\textrm{\'et},*} \mathbb{Z}_{\ell}(1))$, and define $c^{\textrm{\'et}}_{1, S}(\mathcal{L})$ as the image of $c^{\textrm{\'et}}_{1}(\mathcal{L})$ under $g$.

Letting $s \in S$ be a (geometric) point, one defines $r, b$ and $h$ analogously to \autoref{globaldRchern}, and similarly checks that $h = b \circ r \circ g$. 
\end{proof}

\begin{defn}
\label{bigdef}
Let $f : X \to S$ be a smooth projective morphism of $K$-varieties, with $S$ smooth and fibres of dimension $n$. Let $\mathcal{H}^{2k} = R^{2k} f_{*} \Omega^{\bullet}_{X/S}$, $\mathbb{V}^{2k}_{\ell} = R^{2k} f^{\textrm{\'et}}_{*} \mathbb{Z}_{\ell}(k)$ and $\mathbb{V}^{2k} = R^{2k} f^{\textrm{an}}_{*} \mathbb{Z}(k)$. Let $\mathcal{L}$ be a very ample bundle over $S$, and let $\omega^{\textrm{dR}} = c^{\textrm{dR}}_{1, S}(\mathcal{L})$, $\omega^{\textrm{\'et}} = c^{\textrm{\'et}}_{1, S_{\overline{K}}}(\mathcal{L})$ and $\omega$ be the analytification of $\omega^{\textrm{dR}}$. We define: 
\begin{itemize}
\item[(i)] the operators 
\begin{align*}
L^{\textrm{dR}} &: \mathcal{H}^{2k} \to \mathcal{H}^{2k+2} & \hspace{2em} \beta &\mapsto \beta \wedge \omega^{\textrm{dR}} \\
L^{\textrm{\'et}}_{\ell} &: \mathbb{V}^{2k}_{\ell} \to \mathbb{V}^{2k+2}_{\ell} & \hspace{2em} \beta &\mapsto \beta \wedge \omega^{\textrm{\'et}}  \\ 
L &: \mathbb{V}^{2k} \to \mathbb{V}^{2k+2} & \hspace{2em} \beta &\mapsto \beta \wedge \omega ;
\end{align*}
\item[(ii)] and the subbundle and subsystems
\begin{align*}
\mathcal{H}^{2k}_{\textrm{prim}} &= \textrm{ker} \, L^{dR, n-2k+1} \\
\mathbb{V}^{2k}_{\ell, \textrm{prim}} &= \textrm{ker} \, L^{\textrm{\'et}, n-2k+1}_{\ell} \\
\mathbb{V}^{2k}_{\textrm{prim}} &= \textrm{ker} \, L^{n-2k+1} .
\end{align*}
\end{itemize}
Now fix isomorphisms 
\begin{align*}
\xi^{\textrm{dR}} &: R^{2n} f_{*} \Omega^{\bullet}_{X/S} \xrightarrow{\sim} \mathcal{O}_{S} \\
\xi^{\textrm{\'et}}_{\ell} &: R^{2n} f_{\textrm{\'et},*} \mathbb{Z}_{\ell}(n) \xrightarrow{\sim} \mathbb{Z}_{\ell} \\
\xi &: R^{2n} f_{*} \mathbb{Z}(n) \xrightarrow{\sim} \mathbb{Z} ,
\end{align*}
compatible with the comparisons coming from analytification; for instance, using the relative version of the trace isomorphism (see \cite{hartshornedr} and \cite{conradet}). We then further define
\begin{itemize}
\item[(iii)] the polarizations
\begin{align*}
\mathcal{Q}^{2k} &: \mathcal{H}^{2k}_{\textrm{prim}} \otimes \mathcal{H}^{2k}_{\textrm{prim}} \to \mathcal{O}_{S} &  \hspace{2em} \xi^{\textrm{dR}} \circ [ (\alpha, \beta) &\mapsto \alpha \wedge \beta \wedge (\omega^{\textrm{dR}})^{n-2k} ] \\
Q^{2k}_{\ell} &: \mathbb{V}^{2k}_{\ell, \textrm{prim}} \otimes \mathbb{V}^{2k}_{\ell, \textrm{prim}} \to \mathbb{Z}_{\ell} & \xi^{\textrm{\'et}}_{\ell} \circ [ (\alpha, \beta) &\mapsto \alpha \wedge \beta \wedge (\omega^{\textrm{\'et}})^{n-2k} ]  \\ 
Q^{2k} &: \mathbb{V}^{2k}_{\textrm{prim}} \otimes \mathbb{V}^{2k}_{\textrm{prim}} \to \mathbb{Z} & \xi \circ [ (\alpha, \beta) &\mapsto \alpha \cup \beta \cup \omega^{n-2k} ] .
\end{align*}
\end{itemize}
\end{defn}

\begin{prop}
The polarizations $\mathcal{Q}^{2k}, Q^{2k}_{\ell}$ and $Q^{2k}$ of \autoref{bigdef}(iii) satisfy the properties (v), (vi) and (vii) of \autoref{propertiessec}.
\end{prop}

\begin{proof}
For properties (v) and (vi), this follows from the compatibility of the cup product with the comparison isomorphisms, as well as the compatibility of the maps $\xi^{\textrm{dR}}, \xi^{\textrm{\'et}}_{\ell}$ and $\xi$. For property (vii) this follows from the compatibility of the cup product with conjugation by $\textrm{Aut}(\mathbb{C}/K)$, as well as the fact that the section $\omega$ is defined over $K$.
\end{proof}

\subsection{The Absolute Hodge Conjecture}
\label{abshodgesec}

Let $Y$ be a smooth complex projective variety, and $k$ a non-negative integer. We define by $H^{2k}(Y)$ the intersection $H^{2k}(\an{Y}, \mathbb{Q}(k)) \cap H^{k, k}(\an{Y})$, where $H^{2k}(\an{Y}, \mathbb{C}(k)) \cong \bigoplus_{p+q=2k} H^{p, q}(\an{Y})$ is the Hodge decomposition. According to \cite{ChapterElevenNotesonAbsoluteHodgeClasses}, the absolute Hodge conjecture says the following:

\begin{conj}[Absolute Hodge]
\label{abshodge}
Let $\sigma \in \textrm{Aut}(\mathbb{C}/\mathbb{Q})$, and let $c^{\textrm{dR}}_{Y} : H^{2k}(\an{Y}, \mathbb{C}(k)) \xrightarrow{\sim} H^{2k}_{\textrm{dR}}(Y/\mathbb{C})$ and $c^{\ell}_{Y} : H^{2k}(\an{Y}, \mathbb{Q}_{\ell}(k)) \xrightarrow{\sim} H^{2k}_{\textrm{\'et}}(Y, \mathbb{Q}_{\ell}(k))$ be the comparison isomorphisms. Then 
\begin{itemize}
\item[(i)] if $v^{\textrm{dR}} \in H^{2k}_{\textrm{dR}}(Y/\mathbb{C})$ is the image of $v \in H^{2k}(Y)$ under $c^{\textrm{dR}}_{Y}$ then $v^{\textrm{dR}}_{\sigma} \in H^{2k}_{\textrm{dR}}(Y_{\sigma}/\mathbb{C})$ is the image of some $v_{\sigma} \in H^{2k}(Y_{\sigma})$ under $c^{\textrm{dR}}_{Y_{\sigma}}$;
\item[(ii)] if $v^{\ell} \in H^{2k}_{\textrm{\'et}}(Y, \mathbb{Q}_{\ell}(k))$ is the image of $v \in H^{2k}(Y)$ under $c^{\ell}_{Y}$ then $v^{\ell}_{\sigma} \in H^{2k}_{\textrm{\'et}}(Y_{\sigma}, \mathbb{Q}_{\ell}(k))$ is the image of some $v_{\sigma} \in H^{2k}(Y_{\sigma})$ under $c^{\ell}_{Y_{\sigma}}$.
\end{itemize}
\end{conj}

\begin{rem}
Deligne in \cite{Deligne1982} requires the class $v_{\sigma}$ to be \emph{canonical}, in the sense that it is the same for all comparison isomorphisms.
\end{rem}

\begin{rem}
\autoref{abshodge} generalizes to all Hodge tensors in the obvious way.
\end{rem}

\begin{defn}
More generally, we call the statement of \autoref{abshodge} with $\textrm{Aut}(\mathbb{C}/\mathbb{Q})$ replaced by $\textrm{Aut}(\mathbb{C}/K)$ the $K$-absolute Hodge conjecture.
\end{defn}

\begin{defn}
By the ($K$-)\emph{absolute Hodge conjecture for }$v$, we mean the statement of \autoref{abshodge} for some fixed vector $v \in H^{2k}(Y)$.
\end{defn}

Let us now deduce \autoref{mainthm2} from \autoref{mainthm}.

\begin{prop}
\label{actuallyhodge}
In the situation of \autoref{mainthm}, if each element of $W$ is Hodge, then the $K$-absolute Hodge conjecture holds for each element of $W$.
\end{prop}

\begin{proof}
What needs to be checked is that if $w_{\sigma}$ is a rational conjugate to $w \in W$, then $w_{\sigma}$ is automatically Hodge. We recall that a rational vector inside $\mathbb{V}_{s}$ is Hodge if and only if it lies inside $F^{0}$, where $F^{\bullet}$ is the Hodge filtration. As $F^{\bullet}$ is a filtration by $K$-algebraic bundles on $\mathcal{H}$, this shows the result in the de Rham case. In the $\ell$-adic case, we use the fact that $W_{\sigma} = r(W)$ is independent of the point $r$ of $\mathcal{I}_{T}(s, Z, \sigma)$, hence taking $r = r^{\ell}_{\sigma}$ we see that the result holds $\ell$-adically as well.
\end{proof}

\section{Conjugation Isomorphisms and Monodromy}
\label{pointsofIsec}

In this section we establish \autoref{torsorpoints}. The variety $\mathcal{I}_{T}(s, Z, \sigma) \subset \textrm{Isom}(\mathbb{V}_{s}, \mathbb{V}_{s_{\sigma}})$ is defined by the conditions conditions $r \circ \mono_{Z} \circ r^{-1} = \mono_{Z_{\sigma}}$ and $r(t) = t_{\sigma}$ for each $t \in T$, where $r$ is a point of $\textrm{Isom}(\mathbb{V}_{s}, \mathbb{V}_{s_{\sigma}})$ and the first equality is in the sense of $\mathbb{Q}$-subschemes. That $r^{\textrm{dR}}_{\sigma}$ and $r^{\ell}_{\sigma}$ satisfy the second condition is just the assumption that each $t \in T$ has a canonical rational conjugate $t_{\sigma}$, so \autoref{torsorpoints} reduces to the following statement:

\begin{prop}
\label{torsorpoints2}
~
\begin{itemize}
\item[(i)] The map
\[ r^{\textrm{dR}}_{\sigma} : \mathbb{V}_{\mathbb{C}, s} \xrightarrow{c^{\textrm{dR}}_{s}} \mathcal{H}_{\mathbb{C}, s} \xrightarrow{\sigma} \mathcal{H}_{\mathbb{C}, s_{\sigma}} \xrightarrow{(c^{\textrm{dR}}_{s_{\sigma}})^{-1}} \mathbb{V}_{\mathbb{C}, s_{\sigma}} \]
induced by $\sigma$ is satisfies the property that $r^{\textrm{dR}}_{\sigma} \circ \mono_{Z, \mathbb{C}} \circ \left(r^{\textrm{dR}}_{\sigma}\right)^{-1} = \mono_{Z_{\sigma}, \mathbb{C}}$. 
\item[(ii)] For each $\ell$, the map 
\[ r^{\ell}_{\sigma} : \mathbb{V}_{\mathbb{Q}_{\ell}, s} \xrightarrow{c^{\ell}_{s}} \mathbb{V}_{\ell, \mathbb{Q}_{\ell}, s} \xrightarrow{\sigma} \mathbb{V}_{\ell, \mathbb{Q}_{\ell}, s_{\sigma}} \xrightarrow{(c^{\ell}_{s_{\sigma}})^{-1}} \mathbb{V}_{\mathbb{Q}_{\ell}, s_{\sigma}} \]
induced by $\sigma$ satisfies the property that induced by $\sigma$ is satisfies the property that $r^{\ell}_{\sigma} \circ \mono_{Z, \mathbb{Q}_{\ell}} \circ \left(r^{\ell}_{\sigma}\right)^{-1} = \mono_{Z_{\sigma}, \mathbb{Q}_{\ell}}$. 
\end{itemize}
\end{prop}

In both sections that follow, we denote the normalization of an algebraic variety $Z$ by $Z^{n}$. If $Z$ is a subvariety of $S_{\mathbb{C}}$ and $s \in S(\mathbb{C})$ is a point lying in $Z(\mathbb{C})$, we will denote by $s^{n}$ a lift of $s$ to $Z^{n}$.

\subsection{The de Rham Case}
\label{derhamconjproof}

The required statement is implicit in \cite{KOU}; we summarise the argument for expository purposes.

Let us temporarily denote by $\mathbb{V}$ the associated complex local system. We may argue as in the proof of \cite[Proposition 3.1]{KOU} to obtain an equivalence of neutral Tannakian categories
\[ \langle \restr{\mathbb{V}}{Z^{n}} \rangle^{\otimes} \simeq_{\tau} \langle \restr{\mathbb{V}}{Z^{n}_{\sigma}} \rangle^{\otimes} , \]
where the notation $\langle - \rangle^{\otimes}$ denotes the neutral Tannakian category generated by the enclosed object (inside the appropriate category of local systems). Note that in the notation of \cite{KOU} we have $\mathbb{V}^{\sigma} = \mathbb{V}$ from the fact that the connection $\nabla$ is defined over $K$. If $\textrm{Fib}_{Z^{n}, s^{n}}$ and $\textrm{Fib}_{Z^{n}_{\sigma}, s^{n}_{\sigma}}$ are the obvious fibre functors, then the Tannakian groups associated to $\textrm{Fib}_{Z^{n}, s^{n}}$ and $\textrm{Fib}_{Z^{n}_{\sigma}, s^{n}_{\sigma}}$ have natural faithful representations on $\mathbb{V}_{s}$ and $\mathbb{V}_{s_{\sigma}}$, as follows from the fact that an automorphism of $\textrm{Fib}_{Z^{n}, s^{n}}$ (resp. an automorphism of $\textrm{Fib}_{Z_{\sigma}^{n}, s^{n}_{\sigma}}$) is determined by its induced automorphism of $\mathbb{V}_{s^{n}}$ (resp. $\mathbb{V}_{s^{n}_{\sigma}}$). 

Under these representations, the Tannakian groups associated to $\textrm{Fib}_{Z^{n}, s^{n}}$ and $\textrm{Fib}_{Z^{n}_{\sigma}, s^{n}_{\sigma}}$ agree with the complex algebraic monodromy groups of $Z$ and $Z_{\sigma}$. It follows from the fact that $\tau$ is an equivalence that the induced map $r^{\textrm{dR}}_{s} : \mathbb{V}_{\mathbb{C}, s} \xrightarrow{\sim} \mathbb{V}_{\mathbb{C}, s_{\sigma}}$ (the notation is in agreement with our previous definition of $r^{\textrm{dR}}_{s}$) gives an isomorphism of representations, and hence conjugates the (complex) algebraic monodromy groups.

\subsection{The $\ell$-adic Case}
\label{ladicconjproof}

We will reduce the problem to a statement about $\ell$-adic local systems and the \'etale fundamental group, which can then be solved entirely algebraically.

\begin{defn}
Let $\mathbb{V}_{\ell}$ be an $\ell$-adic local system on the complex algebraic variety $Z$. We define $\mono_{\ell, Z}$ to be the identity component of the Zariski closure of $\pi_{1}^{\textrm{\'et}}(Z^{n}, s^{n})$ inside $\textrm{Aut}(\mathbb{V}_{\ell, s})_{\mathbb{Q}_{\ell}}$.
\end{defn}

\begin{lem}
\label{etaleclosuresame}
Let $\mathbb{V}$ be a $\mathbb{Z}$-local system on $\an{Z}$, let $\mathbb{V}_{\ell}$ be an $\ell$-adic local system on $Z$, and suppose that we have an isomorphism $\an{\mathbb{V}}_{\ell} \simeq \mathbb{V} \otimes \mathbb{Z}_{\ell}$. Then the induced isomorphism on fibres identifies the groups $\mono_{Z, \mathbb{Q}_{\ell}}$ and $\mono_{\ell, Z}$.
\end{lem}

\begin{proof}
Using the analytic (resp. $\ell$-adic) equivalence between local systems and monodromy representations, and the canonical identification of $\pi_{1}^{\textrm{\'et}}(Z^{n}, s^{n})$ with the profinite completion $\widehat{\pi_{1}(Z^{n}, s^{n})}$, the statement amounts to the following claim: given the commuting diagram
\begin{center}
\begin{tikzcd}
\pi_{1}(Z^{n}, s^{n}) \arrow[d] \arrow[r] & \textrm{Aut}(\mathbb{V}_{s} \otimes \mathbb{Z}_{\ell}) \arrow[d] \\
\widehat{\pi_{1}(Z^{n}, s^{n})} \arrow[r] & \textrm{Aut}(\mathbb{V}_{\ell, s}) ,
\end{tikzcd} 
\end{center}
the Zariski closures $M$ and $\widehat{M}$ of the images $\textrm{im} (\pi_{1}(Z^{n}, s^{n}))$ and $\textrm{im} (\widehat{\pi_{1}(Z^{n}, s^{n})})$ inside $\textrm{Aut}(\mathbb{V}_{\ell, s})_{\mathbb{Q}_{\ell}}$ coincide. We clearly have $M \subset \widehat{M}$ from the inclusion $\textrm{im} (\pi_{1}(Z^{n}, s^{n})) \subset \textrm{im} (\widehat{\pi_{1}(Z^{n}, s^{n})})$, so it suffices to show the reverse inclusion.

It follows from the fact that a group is dense in its profinite completion that if $\pi_{1}(Z^{n}, s^{n}) \to F$ is a morphism to a finite group $F$, then the induced morphism $\widehat{\pi_{1}(Z^{n}, s^{n})} \to F$ has the same image. As a consequence, the image of $\widehat{\pi_{1}(Z^{n}, s^{n})}$ inside $\textrm{Aut}(\mathbb{V}_{\ell, s}) \cong \varprojlim_{k} \GL(\mathbb{Z}/\ell^{k} \mathbb{Z})$ consists of compatible sequences $(\alpha_{k})_{k \geq 1}$ where each $\alpha_{k}$ is a reduction modulo $\ell^{k}$ of an element in $\textrm{im} (\pi_{1}(Z^{n}, s^{n}))$. Letting $f$ be a polynomial function vanishing on $\textrm{im} (\pi_{1}(Z^{n}, s^{n}))$ with coefficients in $\mathbb{Q}_{\ell}$, it now suffices to show that $f$ vanishes on such compatible sequences $(\alpha_{k})_{k \geq 1}$. Scaling if necessary, we may assume that $f$ has coefficients in $\mathbb{Z}_{\ell}$. But then $f(\alpha_{k}) = 0$ modulo $\ell^{k}$ holds for all $k$, so the result follows.
\end{proof}

Applying \autoref{etaleclosuresame}, the proof of  \autoref{torsorpoints2}(ii) is reduced to showing that the conjugation isomorphism $\mathbb{V}_{\ell, s} \xrightarrow{\sim} \mathbb{V}_{\ell, s_{\sigma}}$ sends the image of $\pi_{1}^{\textrm{\'et}}(Z^{n}, s^{n})$ inside $\textrm{Aut}(\mathbb{V}_{\ell, s})$ to the image of $\pi_{1}^{\textrm{\'et}}(Z^{n}_{\sigma}, s^{n}_{\sigma})$ inside  $\textrm{Aut}(\mathbb{V}_{\ell, s_{\sigma}})$. Applying the isomorphism $\mathbb{V}_{\ell} \simeq \tau_{\sigma}^{*} \mathbb{V}_{\ell}$, it suffices to show that there exists an isomorphism $j : \pi_{1}^{\textrm{\'et}}(S_{\mathbb{C}}, s) \xrightarrow{\sim} \pi_{1}^{\textrm{\'et}}(S_{\mathbb{C}}, s_{\sigma})$ sending $\textrm{im}(\pi_{1}^{\textrm{\'et}}(Z^{n}, s^{n})) \subset \pi_{1}^{\textrm{\'et}}(S_{\mathbb{C}}, s)$ to $\textrm{im}(\pi_{1}^{\textrm{\'et}}(Z^{n}_{\sigma}, s^{n}_{\sigma})) \subset \pi_{1}^{\textrm{\'et}}(S_{\mathbb{C}}, s_{\sigma})$ and a commuting diagram of the following form:
\begin{equation}
\label{fundgpdiag}
\begin{tikzcd}
\pi_{1}^{\textrm{\'et}}(S_{\mathbb{C}}, s) \arrow[d, "j"] \arrow[r] & \textrm{Aut}(\mathbb{V}_{\ell, s}) \arrow[d] \\
\pi_{1}^{\textrm{\'et}}(S_{\mathbb{C}}, s_{\sigma}) \arrow[r] &  \textrm{Aut}((\tau_{\sigma}^{*} \mathbb{V}_{\ell})_{s_{\sigma}}) ,
\end{tikzcd}
\end{equation}
where the horizontal arrows are the natural representations and the vertical arrow on the right comes from the natural isomorphism $\mathbb{V}_{\ell, s} \simeq (\tau_{\sigma}^{*} \mathbb{V}_{\ell})_{s_{\sigma}}$ of \'etale stalks.

\paragraph{Conjugation Isomorphisms:} To describe the map $j$, we begin with a more general construction. We let $X$ be a complex variety, and $\sigma \in \textrm{Aut}(\mathbb{C}/K)$. Then $\sigma$ defines a categorical equivalence $\tau_{\sigma} : \textrm{F\'ET}(X) \to \textrm{F\'ET}(X_{\sigma})$. Denote by $\textrm{Fib}_{X, x}$ the fibre functor at $x \in X$. Given $\alpha \in \pi_{1}^{\textrm{\'et}}(X, x)$, we may define an automorphism $\alpha_{\sigma}$ of $\textrm{Fib}_{X_{\sigma}, x_{\sigma}}$ as follows: for each cover of $X_{\sigma}$, choose an isomorphic cover $X'_{\sigma} \to X_{\sigma}$ in the essential image of $\tau_{\sigma}$; this is the conjugate of a cover $X' \to X$ of $X$. Then define $\alpha_{\sigma}(x'_{\sigma}) = \alpha(x')_{\sigma}$. One checks that $\alpha_{\sigma}$ extends uniquely to a well-defined automorphism $\alpha_{\sigma}$ of $\textrm{Fib}_{X_{\sigma}, s_{\sigma}}$, and that the map $\pi_{1}^{\textrm{\'et}}(X, x) \to \pi_{1}^{\textrm{\'et}}(X_{\sigma}, x_{\sigma})$ defined by $\alpha \mapsto \alpha_{\sigma}$ is a group homomorphism. We define the map $j$ to be the case with $X = S_{\mathbb{C}}$ and $x = s$, and the map $j_{Z^{n}}$ to be the case $X = Z^{n}$ and $x = s^{n}$.


\paragraph{Completing the Proof:} Let $i_{Z^{n}} : \pi_{1}^{\textrm{\'et}}(Z^{n}, s^{n}) \to \pi_{1}^{\textrm{\'et}}(S_{\mathbb{C}}, s)$ and $i_{Z^{n}_{\sigma}} : \pi_{1}^{\textrm{\'et}}(Z^{n}_{\sigma}, s^{n}) \to \pi_{1}^{\textrm{\'et}}(S_{\mathbb{C}}, s)$ be the natural maps. It now suffices check that $j \circ \iota_{Z^{n}} = \iota_{Z^{n}_{\sigma}} \circ j_{Z^{n}}$ and that diagram (\ref{fundgpdiag}) commutes. In the first case, one immediately checks that if $\alpha \in \pi_{1}^{\textrm{\'et}}(S_{\mathbb{C}}, s)$ acts on fibres of $S' \to S_{\mathbb{C}}$ above $s$ by base changing to a cover of $Z^{n}$ and acting via $\pi_{1}^{\textrm{\'et}}(Z^{n}, s^{n})$, then the same is true for $\alpha_{\sigma}$ with respect to $Z^{n}_{\sigma}$ and $\pi_{1}^{\textrm{\'et}}(Z^{n}_{\sigma}, s^{n}_{\sigma})$. In the second case, 
it is immediate from the explicit description of the isomorphism $t : \mathbb{V}_{\ell, s} \xrightarrow{\sim} (\tau_{\sigma}^{*} \mathbb{V}_{\ell})_{s_{\sigma}}$, that together with the map $j$, the map $t$ gives a map of representations (i.e., acting by $\alpha$ then applying $t$ is the same as applying $t$ and then acting by $j(\alpha)$). But $t$ giving a map of representations is equivalent to the commutativity of (\ref{fundgpdiag}).

\section{Variants of the Main Theorem}
\label{mainthmsec}

Although our applications will use \autoref{mainthm} and \autoref{mainthm2}, we wish to briefly explain how \autoref{torsorpoints} may be used to prove variants of \autoref{mainthm} and \autoref{mainthm2} which may be of independent interest. Let us first generalize the notation $\mathcal{I}_{T}(s, Z, \sigma)$ established in the introduction.

\begin{defn}
Let $s \in S(\mathbb{C})$ be a point, $\{ Z_{i} \}_{i \in I}$ a collection of subvarieties of $S_{\mathbb{C}}$ containing $s$, $\sigma \in \textrm{Aut}(\mathbb{C}/K)$ an automorphism, and $T \subset \mathcal{T}(\mathbb{V}_{s})$ a collection of tensors with canonical rational conjugates. Then we define
\[ \mathcal{I} = \mathcal{I}_{T}(s, \{ Z_{i} \}_{i \in I}, \sigma) \subset \textrm{Isom}(\mathbb{V}_{s}, \mathbb{V}_{s_{\sigma}}) \]
by the property that $r \in \mathcal{I}(\mathbb{C})$ if and only if $r \circ \mono_{Z_{i}, \mathbb{C}} \circ r^{-1} = \mono_{\mathbb{Z}_{i, \sigma}, \mathbb{C}}$ for all $i$, and $r$ sends each $t \in T$ to its canonical rational conjugate $t_{\sigma} \in T_{\sigma}$.
\end{defn}

We observe that the normalizer $N$ in the statement of \autoref{mainthm} is equal to $\mathcal{I}_{T}(s, Z, \textrm{id})$, so we may view $\mathcal{I}_{T}(s, \{ Z \}_{i \in I}, \textrm{id})$ as its natural generalization. The following generalization of \autoref{mainthm} is then immediate from \autoref{torsorpoints2}:

\begin{thm}
\label{mainthm3}
Let $\mathbb{V}$ be a motivic variation of Hodge structure on $S$, and let $s \in S(\mathbb{C})$. Suppose that $T \subset \mathcal{T}(\mathbb{V}_{s})$ is a subspace with canonical rational conjugates, that $Z_{i} \subset S_{\mathbb{C}}$ for $i \in I$ is a collection of complex subvarieties containing $s$, and that $W \subset \mathcal{T}(\mathbb{V}_{s})$ is a subspace containing $T$ such that either
\begin{itemize}
\item[(i)] $\mathcal{I}_{T}(s, \{ Z_{i} \}_{i \in I}, \textrm{id})$ fixes $W$ as a $\mathbb{Q}$-subspace; or
\item[(ii)] $\mathcal{I}_{T}(s, \{ Z_{i} \}_{i \in I}, \textrm{id})$ fixes $W$ pointwise.
\end{itemize}
Then 
\begin{itemize}
\item[(i)] each element $w$ of $W$ has rational conjugates; or
\item[(ii)] each element $w$ of $W$ has a canonical rational conjugate,
\end{itemize}
in cases (i) and (ii) respectively.
\end{thm}

\begin{proof}
As in the proof of \autoref{mainthm}, we obtain from \autoref{torsorpoints2} that $r^{\textrm{dR}}_{\sigma}$ and $r^{\ell}_{\sigma}$ for all $\ell$ are points of $\mathcal{I} = \mathcal{I}_{T}(s, \{ Z_{i} \}_{i \in I}, \sigma)$. It follows as in \autoref{mainthm} that the points of $W$ are defined over $\overline{\mathbb{Q}}$ and (under an appropriate isomorphism $\mathbb{C} \cong \overline{\mathbb{Q}_{\ell}}$) over $\overline{\mathbb{Q}}_{\ell}$ for every $\ell$. By \autoref{affinespptclaim} we conclude in that $r(W_{\sigma})$ is a $\mathbb{Q}$-subspace, where $r \in \mathcal{I}(\mathbb{C})$. In Case (ii) we additionally know that $r(w) = w_{\sigma}$ is independent of the choice of $r$, making the conjugates canonical.
\end{proof}

Finally let us give a similar, but simpler argument which establishes cases of a kind of ``$\overline{\mathbb{Q}}$-absolute Hodge'' conjecture. We note that a ``$\overline{\mathbb{Q}}$-absolute Hodge'' conjecture may in fact be enough for certain applications of absolute Hodge to the algebraicity of periods. The following requires no $\ell$-adic input.

\begin{thm}
\label{mainthm4}
Let $\mathbb{V}$ be a variation of Hodge structure on $S$, suppose that $S$, $\nabla$ and $\mathcal{H}$ are defined over $K$, and choose $s \in S(\mathbb{C})$. Let $T \subset \mathcal{T}(\mathbb{V}_{\overline{\mathbb{Q}}, s})$ be a $\overline{\mathbb{Q}}$-subspace with canonical rational conjugates, let $Z_{i} \subset S_{\mathbb{C}}$ for $i \in I$ be a collection of complex subvarieties containing $s$, and suppose that $W \subset \mathcal{T}(\mathbb{V}_{\overline{\mathbb{Q}}, s})$ is a subspace containing $T$ such that the orbit of $[W]$ under $\mathcal{I} = \mathcal{I}_{T}(s, \{ Z_{i} \}_{i \in I}, \textrm{id})$ inside the appropriate Grassmanian is finite. Then $W$ has $\overline{\mathbb{Q}}$-rational conjugates.
\end{thm}

\begin{proof}
The proof follows immediately, as \autoref{torsorpoints2} ensures that $r^{\textrm{dR}}_{\sigma}$ (and $r^{\ell}_{\sigma}$, if an appropriate comparison to an $\ell$-adic local system exists) lies inside $\mathcal{I}$ and the assumptions ensure that $r(W)$ is defined over $\overline{\mathbb{Q}}$ for any point $r$ of $\mathcal{I}$.
\end{proof}

Finally, we note that there is an additional step involved in translating \autoref{mainthm3} and \autoref{mainthm4} to absolute Hodge statements (like \autoref{mainthm2}), but this follows exactly as in \autoref{abshodgesec} and we leave this to the reader.

\section{Applications to Absolute Hodge}

\subsection{Voisin's Example}
\label{voisinexamplesec}

In this section we prove \autoref{voisinexample}, following an approach laid out by Voisin in Section 3 of \cite{voisin_2007}. Let us first revisit the proof of our \autoref{maincor} in the context of \autoref{mainthm2}.

We let $\mathbb{V}$ be a motivic variation on the smooth $K$-variety $S$, and $T \subset \mathcal{T}(\mathbb{V}_{s})$ be the $\mathbb{Q}$-span of the polarization $Q_{s}$; by \autoref{conjpreservespol} the subspace $T$ has canonical rational conjugates. Suppose that there exists a subvariety $Z \subset S_{\mathbb{C}}$ containing $s$ such that the fixed locus of $\mono_{Z}$  inside $\mathbb{V}^{m,n}_{s}$ is the line spanned by the Hodge class $w \in \mathbb{V}^{m,n}_{s}$. Then if $n$ is a complex point of the normalizer of $\mono_{Z}$ inside $\mathbf{G}_{T}$, we have that $n w = \lambda w$ for some $\lambda \in \mathbb{C}$, and from the fact that $Q(n w, n w) = Q(w, w)$ we learn that $\lambda = \pm 1$. \autoref{mainthm2} therefore applies, and the tensor $w$ is absolutely Hodge.

A very similar argument is given by Voisin in her paper (see Remark 1.2, Theorem 0.5(1) in \cite{voisin_2007}). The language used is slightly different: Voisin only considers Hodge vectors $v$ (the $(m, n) = (0, 1)$ case); the variety $Z$ is taken to be a special subvariety (irreducible component of the Hodge locus); and the condition that the fixed locus of $\mono_{Z}$ consist of exactly the line spanned by $w$ takes the form of the condition that the restricted variation $\restr{\mathbb{V}}{Z}$ has $\mathbb{Q} w$ as its constant subvariation. Voisin also only obtains the weaker conclusion that the conjugate $v_{\sigma}$ is a $\overline{\mathbb{Q}}$-scalar multiple of a Hodge class; the essential difference, it seems to us, is the presence of $\ell$-adic input in our argument.

To study the example of \autoref{voisinexample}, Voisin states a criterion \cite[Theorem 3.1]{voisin_2007}, whose proof gives the following:

\begin{prop}[Voisin]
\label{voisincriterion}
Let $\mathbb{V}$ be a polarizable variation of Hodge structure of weight $2k$. Denote by $\overline{\nabla}^{p,q}$ the $\mathcal{O}_{S}$-linear map $\mathcal{H}^{p, q} \to \mathcal{H}^{p-1,q+1} \otimes \Omega_{S}$ induced by the connection, let $\alpha \in \mathbb{V}_{s}$ be a Hodge class with Hodge locus $Z_{\alpha}$, and let $\lambda$ be its projection to $\mathcal{H}^{k,k}$. Suppose that
\begin{itemize}
\item[(i)] the map $\mu_{\lambda}(v) = \overline{\nabla}_{v}(\lambda)$ is surjective;
\item[(ii)] for $p > k, p+q = 2k$ the restriction of $\overline{\nabla}^{p,q}$ to the the tangent space $T_{s} Z_{\alpha}$ is injective;
\item[(iii)] and the restriction of $\overline{\nabla}^{k,k}$ to the tangent space $T_{s} Z_{\alpha}$ has kernel equal to the span of $\lambda$.
\end{itemize}
Then the fixed locus in $\mathbb{V}_{s}$ of $\mono_{Z_{\alpha}}$ is exactly the line spanned by $\alpha$.
\end{prop}

\begin{proof}[Proof of \autoref{voisinexample}:]
We observe that the variation $\mathbb{V}_{-}$ is motivic: taking anti-invariants defines an appropriate subbundle $\mathcal{H}_{-} \subset \mathcal{H}$ and subsystem $\mathbb{V}_{\ell,-} \subset \mathbb{V}_{\ell}$, compatibly with the comparison isomorphisms, and the required properties are simply obtained by restriction. Arguing as Voisin does in \cite[Example 3.4]{voisin_2007}, we verify that the hypotheses of \autoref{voisincriterion} hold for classes $\alpha$ whose projections $\lambda$ lie in a certain topologically dense subset of the underlying real subbundle of $\mathcal{H}^{k,k}_{-}$. The result then follows from \autoref{mainthm2}.
\end{proof}

\subsection{Absolutely Hodge Isomorphisms}
\label{genisosec}

In this section we study the question of when an isomorphism of Hodge structures between two fibres in a motivic family is absolute Hodge. Related questions are considered in \cite[Section~3]{voisin_2007}. 

\begin{defn}
\label{boxvardef}
Let $S_{1}$ and $S_{2}$ be algebraic varieties, and let $\mathbb{V}_{1}$ and $\mathbb{V}_{2}$ be variations of Hodge structure on $S_{1}$ and $S_{2}$, respectively. We define $\mathbb{V}_{1} \boxtimes \mathbb{V}_{2}$ to be the variation $p^{*}_{1} \mathbb{V}_{1} \otimes p^{*}_{2} \mathbb{V}_{2}$ on $S_{1} \times S_{2}$, where $p_{i} : S_{1} \times S_{2} \to S_{i}$ is the projection.
\end{defn}

Suppose that $\mathbb{V}$ is an integral variation of Hodge structure on $S$, and let $Z \subset S \times S$ be the algebraic subvariety whose points are pairs $(x,y)$ such that $\mathbb{V}_{x}$ is isomorphic to $\mathbb{V}_{y}$ as an integral polarized Hodge structure. These isomorphisms are Hodge tensors of the variation $\mathbb{V} \boxtimes \mathbb{V}^{*}$. When a Torelli theorem is avaliable, the locus $Z$ is simply equal to the diagonal $\Delta$, and its algebraic monodromy is easily determined to be equal to the image of $\mono_{S}$ under the diagonal action. In the general case, let $Z_{0}$ be an irreducible component of $Z$ containing the diagonal $\Delta \subset S \times S$. We determine $\mono_{Z_{0}}$ by reducing to the Torelli case using a recent result of Bakker, Brunebarbe and Tsimerman \cite{OMINGAGA}.

To show our main result we first prove some lemmas on algebraic monodromy:

\begin{lem}
\label{properalgmonolem}
Let $f : S \to T$ be a proper, dominant morphism of irreducible complex algebraic varieties, and let $\mathbb{V}$ be a local system on $T$. Then for any point $s \in S$, the isomorphism $(f^{*} \mathbb{V})_{s} \xrightarrow{\sim} \mathbb{V}_{f(s)}$ induces an isomorphism between $\mono_{S}$ and $\mono_{T}$.
\end{lem}

\begin{proof}
That $f : S \to T$ is proper means that it has a Stein factorization $f = g \circ h$, where $h : S \to V$ is proper with connected fibres and $g$ is finite. We are thus reduced to showing the theorem in two cases:
\begin{itemize}
\item[(i)] when $f$ additionally has connected fibres;
\item[(ii)] and when $f$ is additionally finite.
\end{itemize}

Let us first assume that both $S$ and $T$ are normal. If $X$ is a normal variety, then the fundamental group of any open subvariety $U$ surjects onto the fundamental group of $X$, so we may replace both $S$ and $T$ by Zariski open sets, restricting $f$ appropriately. In case (i) this lets us assume that $f$ is surjective and flat, and hence by \cite[\href{https://stacks.math.columbia.edu/tag/01UA}{Lemma 01UA}]{stacks-project} universally open. Koll\'ar has shown in \cite{pathlifting} that a universally open, surjective morphism of complex varieties with connected fibres satisfies the \emph{two-point path lifting property} (see \cite[Definition 30]{pathlifting}), and therefore induces a surjection of fundamental groups, from which the result follows. In case (ii) we may assume that $f$ is \'etale, from which the result follows as the monodromy of $S$ will be finite index in the monodromy of $T$, and hence have the same Zariski closure.

Working now in the general case, it follows from the definition of algebraic monodromy that the algebraic monodromy of a local system $\mathbb{V}$ on $X$ agrees with the algebraic monodromy of the restriction of $\mathbb{V}$ to the normal locus of $X$. The result then follows by restricting $f$.
\end{proof}

\begin{lem}
\label{monoisdiag}
Let $\mathbb{V}$ be a polarizable integral variation of Hodge structure on the smooth complex algebraic variety $S$, and let $Z_{0}$ and $\Delta$ be as above. Then $\mono_{Z_{0}} = \mono_{\Delta}$.
\end{lem}

\begin{proof}
The statement is unchanged under replacing $S$ with a finite \'etale cover, so we may assume that $S$ has unipotent monodromy at infinity. Let $\varphi : S \to \Gamma \backslash D$ be the period map on $S$; here $D$ is the full period domain of polarized integral Hodge structures on a fixed lattice $V$, and $\Gamma = G(\mathbb{Z})$ where $G = \textrm{Aut}(V, Q)$. Then $\Gamma \backslash D$ may be interpreted as the moduli space for polarized integral Hodge structures of the same type as the fibres of $\mathbb{V}$, and the map $\varphi$ sends $s$ to the isomorphism class of $\mathbb{V}_{s}$. Arguing as in \cite[Corollary 13.7.6]{CMS} we may complete the variety $S$ by adding finitely many points so that $\varphi$ is proper. Applying the main theorem of \cite{OMINGAGA}, we find that there exists a factorization $\varphi = \varphi' \circ g$ where $g : S \to T$ is a dominant proper map of algebraic varieties and $\varphi : T \hookrightarrow \Gamma \backslash D$ is a closed embedding.

The variety $T$ admits a variation $\mathbb{V}'$ such that $\mathbb{V}$ may be identified with $g^{*} \mathbb{V}'$, and $\varphi'$ is the period map of $\mathbb{V}'$. The variety $Z$ is then the inverse image of the diagonal under the product morphism $g \times g : S \times S \to T \times T$. If $\Delta_{S} \subset S \times S$ and $\Delta_{T} \subset T \times T$ are the respective diagonals, the maps $Z_{0} \to \Delta_{T}$ and $\Delta_{S} \to \Delta_{T}$ are both dominant and proper, so the result follows by \autoref{properalgmonolem}.
\end{proof}

Let us note that it is a consequence of Deligne's Principle A (see \cite[Theorem 3.8]{deligne}) that if the tensor $Q_{s}$ has canonical rational conjugates, then so does every rational tensor fixed by $\textrm{Aut}(\mathbb{V}_{s}, Q_{s})$. \autoref{genisothm} now follows from the more general:

\begin{thm}
If $\mathbb{V}$ is a motivic variation of Hodge structure on the smooth $K$-variety $S$ such that rational tensors fixed by $\mono_{S}$ have canonical rational conjugates and the points of the center $\mathcal{Z}(\mono_{S})$ are defined over $\mathbb{Q}$, then tensors fixed by $\mono_{Z_{0}}$ have rational conjugates, and generic isomorphisms between $\mathbb{V}_{x}$ and $\mathbb{V}_{y}$ are $K$-absolutely Hodge.
\end{thm}

\begin{proof}Let $T \subset \mathcal{T}(\mathbb{V}_{s})$ be the subspace of $\mathbb{Q}$-tensors fixed by $\mono_{S}$, and let $T^{*}$ be its dual. As the subspace $T$ (resp. the subspace $T^{*}$) is fixed under algebraic monodromy, parallel transport gives a path-independent translate of $T$ (resp. $T^{*}$) to any fibre of $\mathbb{V}$ (resp. $\mathbb{V}^{*}$). It follows that we may unambiguously refer to the subspace $R = T \otimes T^{*}$ inside $\mathcal{T}(\mathbb{V}_{x} \otimes \mathbb{V}^{*}_{y})$ for any pair $(x, y)$. The fact that $T$ has canonical rational conjugates implies the same fact for $R$. We then have:

\begin{claim}
The groups $\mathbf{G}_{R}$ and $\mathbf{G}_{T} \times \mathbf{G}_{T^{*}}$ are equal as subgroups of $\textrm{Aut}(\mathbb{V}_{x} \otimes \mathbb{V}^{*}_{y})$.
\end{claim}

\begin{proof}
Let $U \subset \mathcal{T}(\mathbb{V}_{x} \otimes \mathbb{V}^{*}_{y})$ be the subspace of tensors fixed by $\mathbf{G}_{T} \times \mathbf{G}_{T^{*}}$. It suffices to show that $U = T \otimes T^{*}$. Let us first argue that the space $U$ is spanned by pure tensors; i.e., it has a basis of the form $u \otimes u'$ where $u \in \mathcal{T}(\mathbb{V}_{x})$ and $u' \in \mathcal{T}(\mathbb{V}^{*}_{y})$. Indeed, suppose that $\mathbf{G}_{T} \times \mathbf{G}_{T^{*}}$ fixes $\sum_{k = 1}^{m} v_{k} \otimes v'_{k}$, where $\{ v_{k}  \}_{k = 1}^{m}$ and $\{ v'_{k} \}_{k = 1}^{m}$ form linearly independent sets. Let $(g, g') \in \mathbf{G}_{T}(\mathbb{C}) \times \mathbf{G}_{T^{*}}(\mathbb{C})$ be a point. Then from the fact that $(g, 1)$ fixes $\sum_{k = 1}^{m} v_{k} \otimes v'_{k}$ we learn that $g$ fixes each $v_{k}$, and analogously we learn that $(1, g')$ fixes each $v'_{k}$. It follows that $(g, g')$ fixes each term in the sum.

We are reduced to arguing that any element of the form $v \otimes v'$ with $v \in \mathcal{T}(\mathbb{V}_{x})$ and $v' \in \mathcal{T}(\mathbb{V}^{*}_{y})$ fixed by $\mathbf{G}_{T} \times \mathbf{G}_{T^{*}}$ must satisfy $v \in T$ and $v' \in T^{*}$. From the fact that $\mathbf{G}_{T} \times \mathbf{G}_{T^{*}}$ acts on $v$ and $v'$ separately we learn immediately that $\mathbf{G}_{T}$ (resp. $\mathbf{G}_{T^{*}}$) acts on $v$ (resp. $v'$) through a character. As $\mathbf{G}_{T}$  (resp. $\mathbf{G}_{T^{*}}$) is a connected, semisimple group (isomorphic to $\mono_{S}$), we learn that $\mathbf{G}_{T}$ (resp. $\mathbf{G}_{T^{*}}$) fixes $v$ (resp. $v'$).
\end{proof}

We may now complete the proof by computing the normalizer $N$ of $\mono_{Z_{0}}$ inside $\mono_{S} \times \mono_{S}$. It follows from \autoref{monoisdiag} that $\mono_{Z_{0}}$ is the diagonal subgroup $D$ of $\mono_{S} \times \mono_{S}$. The normalizer $N$ of $D$ is the group generated by $D$ and the product of centers $\mathcal{Z}(\mono_{S}) \times \mathcal{Z}(\mono_{S})$. It follows that $N$ preserves the subspace $W$ of all rational tensors fixed by $\mono_{Z_{0}}$, and so the result follows by \autoref{mainthm} and \autoref{mainthm2}.
\end{proof}

\bibliographystyle{alpha}
\bibliography{hodgeloci}

\end{document}